\documentclass{elsart}

 \usepackage{graphicx}

\usepackage{amssymb}

\usepackage[T1]{fontenc}
\usepackage[latin1]{inputenc}

\usepackage[english,francais]{babel}
\usepackage{psfrag}
\usepackage{verbatim}

\newtheorem{e-proposition}[theorem]{Proposition}

\newtheorem{e-definition}[theorem]{Definition\rm}

\def\og{\leavevmode\raise.3ex\hbox{$\scriptscriptstyle\langle\!\langle$~}}
\def\fg{\leavevmode\raise.3ex\hbox{~$\!\scriptscriptstyle\,\rangle\!\rangle$}}

\newcommand{\am}             [0]     {a_{max}}

\newcommand{\bj}             [0]     {{\bf J}}

\newcommand{\bl}             [0]     {{\bf L}}
\newcommand{\bq}             [0]     {{\bf Q}}
\newcommand{\bu}             [0]     {{\bf U}}

\newcommand{\br}             [0]     {{\bf R}}

\begin{document}

\begin{frontmatter}


\selectlanguage{english}
\vspace*{-0.75cm}
\title{\hbox{A high order purely frequency-based harmonic balance} \hbox{formulation for continuation of periodic solutions}}

\vspace*{-2.25cm}
\author[LMA]{Bruno Cochelin}
\ead{bruno.cochelin@ec-marseille.fr}
\author[LMA]{Christophe Vergez\corauthref{cor1}}
\ead{vergez@lma.cnrs-mrs.fr}
\corauth[cor1]{Corresponding author. Tel.: +33 491 16 41 63 ; Fax: +33 491 16 40 12}

\vspace*{-0.25cm}
\address[LMA]{Laboratoire de Mécanique et d'Acoustique (LMA-CNRS),\\ 31 Chemin Joseph Aiguier, 13402 Marseille Cedex 20, France}

\vspace*{-0.75cm}
\begin{abstract}
Combinig the harmonic balance method (HBM) and a continuation method is a well-known technique  to follow 
the periodic solutions of dynamical systems when a control parameter is varied. However, since deriving
 the algebraic system containing the Fourier coefficients can be a highly cumbersome procedure, the classical HBM is often
limited to polynomial (quadratic and  cubic) nonlinearities  and/or a few harmonics. Several variations on  the classical HBM, such as the incremental HBM or the alternating frequency/time domain HBM,  have been presented in the literature to overcome this shortcoming. Here, we present
an alternative approach that can be applied to a very large class of dynamical systems (autonomous or forced) with smooth equations. The main idea is to systematically recast the dynamical system in quadratic polynomial form before applying the HBM. Once the equations have been rendered  quadratic,
it becomes obvious to derive  the algebraic system and solve it by the so-called ANM continuation technique.
Several classical examples are presented to illustrate the use of this numerical approach. 
\end{abstract}

\end{frontmatter}

\selectlanguage{english}
\section{Introduction}

In the field of engineering, especially when dealing with nonlinear vibrations, it is often required to
compute the periodic solutions of a system of nonlinear differential equations (\cite{nayfeh79,szemplinska90}). 
In the following, we focus on numerical methods to find periodic solutions. They are usually subdivided into two main groups : those relying on the time-domain formulation and those relying on the frequency-domain formulation. The first group  consists of  methods where a time integration algorithm, which is generally limited to a single period, is used to transform the original differential system into a system of algebraic equations, which are then solved by continuation. With these methods, the unknowns in the algebraic system are the values of the  original unknown variables at grid points along the periodic orbit. The classical shooting technique (\cite{Seydel94,nayfeh95}) and the orthogonal collocation methods used in AUTO (\cite{Doedel91,Doedel07}) belong to this group. The second group corresponds to the so-called harmonic balance method where the unknown variables are decomposed into truncated Fourier series. In this case, the unknowns in the final algebraic system, which is obtained by balancing the harmonics featuring in the differential equations, are simply the Fourier coefficients of the original unknowns. Methods of both kinds are widely used in many applications and research is still required to improve them. In  practical terms, the choice between the time-domain or the frequency-domain approach depends mainly on whether the periodic solution can be decomposed with a few Fourier components. The method chosen to deal with a given problem can therefore depend on the operating point.

In this paper, we focus on the harmonic balance method  (HBM) and some of its variations, which are briefly recalled here. Although the name "harmonic balance" seems to have first appeared in 1936 (\cite{krylov36}), this method has been widely used only since the sixties, and especially for electrical and mechanical engineering purposes. Forced vibrations were first studied (\cite{urabe65}), and self-sustained oscillations a decade later (\cite{stokes72}). The study by  Nakhla \& Vlach (\cite{nakhla76}) is often said to be a milestone in the modern formulation of HBM. 
The classical HBM is simple in its principle, but it can be cumbersome or even impracticable,  as stated by Peng et al. (\cite{peng2008}) when
the system contains complex nonlinearities and when a large number of harmonics is required, ie, more than 5 or 10.
The first problem which arises  is how to derive  the algebraic system for the Fourier coefficients and the second problem is how  to solve 
 this strongly nonlinear system efficiently. To overcome these shortcomings, many variations on the basic HBM have been proposed in the literature. Here, we will simply give an overall picture of  incremental harmonic balance (IHB) and  alternating frequency/time domain harmonic balance (AFT) method. With the IHB  method (\cite{lau1981}), the incremental-iterative method used for the continuation is closely combined with the harmonic principle, so as to be able to apply the HB principle to the incremental linear problem instead of the nonlinear original system.  Nonlinearity of all kinds can be treated using this technique \cite{iu1986,pierre1985}.  In the AFT variant, the harmonic balance of the Fourier coefficients is not  explicitly performed.
At each increment (or iteration) in the continuation procedure, the unknowns are transferred to the time domain (using the inverse Fast Fourier Transform,which is  denoted IFFT below), so as be able to use the original system of equations. The nonlinear responses are then transferred back to the frequency domain using the Fast Fourier Transform, (denoted FFT below). The use of FFT/IFFT procedures seems to have started in \cite{ling1987}. This technique has been used, for instance, to analyse nonlinear vibrations in mechanical systems with contact and dry friction  \cite{navicet2003}. Since the pioneering work on IHB and AFT, many variations have been proposed to extend the applications, as well as to decrease the computational cost (see \cite{lay2008} or \cite{liu2006} for some recent examples).

An alternative strategy is proposed here for applying the classical HBM with a large number of harmonics,
 without avoiding explicit balancing. 
The  idea is to apply a transformation to the original system of differential equations before applying the HBM. The aim of this procedure is to transform the non-linearities present in the original system into purely polynomial quadratic terms. The classical HBM procedure can then be easily applied. This transformation might seem to be a limitation of the method, since not every system can be recast in polynomial quadratic form. However, it will be shown in this paper that a very large class of systems with smooth equations can indeed be recast in quadratic form by making a few  algebraic manipulations and a a few  additions of equations and auxiliary variables. The idea of using a quadratic formulation to simplify the (Fourier) series expansion was inspired by  another numerical method, the so-called  Asymptotic Numerical Method (ANM) (\cite{azrar93,cochelin94a}) which was presented for the first time in 1990 (\cite{damil90}). ANM  is a continuation technique (\cite{cochelin94b}) involving high order powers series expansion of the  branches of solutions.  For reasons that will be given later on, this ANM continuation is the ideal means of  solving the algebraic system resulting  from the HB method proposed here.

The paper is organised as follows:  details of the transformation into a polynomial quadratic form, the  application of the HBM method and the principle of continuation are given in section 2 in the cases of autonomous systems. Three simple examples are given as an illustration. Section 3 is devoted to  periodically forced systems, and two further examples
are presented. Section 4 gives  the results of all the examples and ends with some concluding comments and perspectives.

\section{Presentation of the method for obtaining periodic solutions of autonomous systems }

Let us consider an autonomous  system of  differential equations 
\begin{equation}\label{systdyna}
\dot Y=f(Y,\lambda)
\end{equation}
where $Y$  is a vector of unknowns, $f$ a smooth nonlinear vector valued function and $\lambda$ a real parameter. The dot stands for the 
derivative with respect to time $t$. We assume that this system
 has  branches of periodic solutions when $\lambda$ varies, and we want to find and follow them
by applying the harmonic balance method and a continuation procedure (path following technique). The important case of a forced  (non-autonomous) system will be dealt with separately in section \ref{s:forced}.

\subsection{The harmonic balance principle}

As recalled in the Introduction, the HBM  consists basically in decomposing $Y(t)$ into a truncated Fourier series :
\begin{equation}\label{fourierser}
Y(t)= Y_0 + \sum_{k=1}^H Y_{c,k} \cos(k \omega t) + \sum_{k=1}^H Y_{s,k} \sin(k \omega t)
\end{equation}
This  ansatz is put into Eq. (\ref{systdyna}) and $f(Y,\lambda)$ is expanded into 
Fourier series. By balancing the first $2H+1$ harmonic terms, one obtains an algebraic system
with $2H+1$ vector equations  for the $2H+1$ vector unknowns $Y_i$, the unknown pulsation $\omega$, and the parameter $\lambda$. Adding a phase condition (\cite{Seydel94}) yields a system with $N$ equations and $N+1$ variables.
The branches of solutions of this algebraic system are then followed using a continuation technique.
This procedure provides only approximate periodic solutions, since in the expansion of Eq. (\ref{systdyna}), all the harmonic terms greater than $H$ remain unbalanced. However, if the number of harmonics $H$ is large enough, accurate solutions can be obtained.

The most crucial point is the expansion of the vector $f(Y,\lambda)$ into  Fourier series.
This can be quite a straighforward procedure if $f(Y,\lambda)$ is a polynomial of  degree two or three, and if the number of harmonics $H$ is sufficiently small. In this case, the expansion can be carried out  by hand.  But in the most  general situation, where $f(Y,\lambda)$ shows nonlinearities of any kinds, this computation can be very cumbersome, even with the help of a symbolic software program. The only alternative is then to use numerical procedures to estimate the Fourier series of $f(Y,\lambda)$, by performing iterative FFT/IFFT steps, for example. This drawback was the starting-point in the developpment of many variants of the HBM as mentioned in the Introduction.

\subsection{A key point: the quadratic recast \label{ss:quadratic_recast}}

The main aim of this paper is to present a simple but powerful procedure that overcomes the drawback  mentioned above. The idea is to recast the original system (\ref{systdyna}) into a new system where the nonlinearities are at most quadratic polynomials . The application of the harmonic balance method subsequently becomes quite straighforward.
This  new quadratic system, which  will be written as follows 
\begin{equation}\label{equaquad}
m(\dot Z) = c + l(Z) + q(Z,Z)
\end{equation}
can contain both differential and algebraic equations. $N_e$ is taken to denote the number of equations. The unknown vector $Z$ (size $N_e$) contains the original components of the vector $Y$ and some new
 variables which are  added  to get the quadratic form. The right hand side of Eq. (\ref{equaquad})  is written as follows: $c$ is a constant vector with respect to the unknown $Z$,  $l(.)$ is a  linear vector valued operator with respect to the vector entry, and $q(.,.)$ is a quadratic vector valued operator, which is linear with respect to both entries. At this stage, the three vectors $c$, $l(.)$
 and $q(.,.)$ may depend on the real parameter $\lambda$. However, it is preferable for the quadratic operator $q(.,.)$  not to depend on $\lambda$, which may require  introducting a new variable in $Z$ (see example 1 below). on the left hand side, $m(.)$  is  a linear vector valued operator with respect to the vector entry. The algebraic equations correspond to zero values of $m(.)$.

Since the nonlinearities are quadratic, the HBM can obviously be easily and systematically applied on Eq. (\ref{equaquad}), even with a large number of harmonics. The question is now: can any system be easily put into the form of  Eq. (\ref{equaquad}) ? As we will see, the answer is yes, for a very large class of nonlinear systems.

The idea  of recasting a nonlinear system into a new one with a quadratic polynomial nonlinearity is frequently used when considering  the so-called ANM continuation method, which consists in computing power series expansion of solution branches, see (\cite{azrar93,cochelin94a,cochelin94b,cochelin94c}) for the basic algorithms and (\cite{zahrouni99,azrar99,cadou01}) for some applications in mechanics. The quadratic recast  
gives a very simple and systematic algebra for the series determination, even with high orders of truncature. Similar benefits are obtained with the Fourier series expansions. 

It is now worth illustrating this main idea by giving some elementary examples where the recast from Eq. (\ref{systdyna}) to Eq. (\ref{equaquad}) will be performed explicitely. We will deal first with the classical Van der Pol oscillator, the Rössler system and a model for clarinet-like musical instruments. Other examples, with  nonlinearities of others kinds such as rational fraction, will be given later on.

\medskip

\textbf{Example 1: } The Van der Pol oscillator.
We take the following second order autonomous system, known as the Van der Pol oscillator, where the parameter $\lambda$ governs the amplitude of the nonlinear damping term.
\begin{equation}\label{vdp1}
\ddot u - \lambda(1-u^2) \; \dot u + u = 0
\end{equation}
This equation can be classically recast into  first order system ( as Eq. (\ref{systdyna})) by introducing the
velocity $v(t)= \dot u(t)$ as a new unknown. We obtain
\begin{equation}\label{vdp2}
\begin{array}{ccc}
\dot u & = & v   \\
\dot v & = & \lambda(1-u^2) \; v - u 
\end{array}
\end{equation}
If we now introduce the auxiliary variables $w(t)= 1- u^2(t)$ and
$r(t)= v(t) w(t)$, the system can be written in the form of Eq. (\ref{equaquad}) with  $Z=[ u, v, w, r ]^t$. The various terms are arranged below so that it will be clear to the readers how the different functions $m$, $c$, $l$ and $q$ are formed.  
\begin{equation}\label{vdp3}
\begin{array}{ccccc}
\dot u & = & 0 & +v &  \\
\dot v & =& 0  & -u + \lambda r & \\
0 & =& 1 & -w &- u^2  \\
 \underbrace{0 }_{m(\dot Z)} & = & \underbrace{ \quad \ 0  \quad  }_{c} & \underbrace{ \quad \quad +r  \quad \quad }_{l(Z,\lambda)}  \  & \underbrace{\quad  - v w   \quad }_{q(Z,Z)}  \\
\end{array}
\end{equation}
We finally obtain two first order differential equations and two algebraic equations with quadratic polynomial nonlinearities. Only the operator $l$ depends here on $\lambda$.

 \medskip

\textbf{Example 2: }The Rössler system. We take the following first order autonomous system of three equations known as the Rössler system. 
\begin{equation}\label{rossler1}
\begin{array}{lcl}
\dot x & = & -y -z \\
\dot y & = & x + a y \\
\dot z & = & b + z(x-\lambda)
\end{array}
\end{equation}
where  $x$, $y$, $z$ are functions of time and $a$, $b$, $c$ are float parameters.
This system can  be written in the form of  Eq. (\ref{equaquad}) without any auxiliary variables (i.e. $Z=[x,y,z]^t$):
\begin{equation}\label{rossler2}
\begin{array}{ccccc}
\dot x & = & 0 & -y -z &  \\
\dot y & =& 0  & +x +ay & \\
 \underbrace{\dot z}_{m(\dot Z)} & =& \underbrace{\quad b \quad}_{c}  & \underbrace{\quad \quad -\lambda z \quad \quad}_{l(Z,\lambda)} &\underbrace{\quad + z x   \quad }_{q(Z,Z)}   \\
\end{array}
\end{equation}
 Here again,  only the operator $l$ depends on $\lambda$. This example will be used here to show that the approach presented in this paper can be used to find period-doubling and  bifurcations with period 4.

\medskip

\textbf{Example 3:} A model for clarinet-like musical instruments. In the case of small amplitude oscillations, a simple model for reed instruments (clarinet, saxophone, etc) can be written using a modal formulation (see \cite{Silva_2007}):
\begin{equation}\label{clari1}
\begin{array}{l}
\ddot x + q_r\omega_r \dot x + \omega_r^2 x = \omega_r^2 p \\
\ddot p_n + 2\alpha_n c \, \dot p_n + \omega_n^2 p_n = \frac{2c}{l} \, \dot u \qquad \forall n \in [1 \dots N]\\
u = \zeta(1-\lambda+x)\sqrt{\lambda - p} \\
p = \sum_{n=1}^{N}p_n
\end{array}
\end{equation}
 where the unknowns $x$, $p_{n=1..N}$, $p$, $u$ are functions of time and $q_r$, $\omega_r$, $\alpha_n$, $\omega_n$, $c$, $l$, $\zeta$ are given parameters describing either the physics or the player's action. $N$ is the number of acoustic modes.  $\lambda$ stands here for the blowing pressure. 

It is worth noting that even if this system contains a square root, it can be recast in the form of Eq. (\ref{equaquad}):  if we introduce the auxiliary variables $y=\dot x$, $z_n = \dot p_n$ and $v=\sqrt{\lambda-p}$, system (\ref{clari1}) can be rewritten, with $Z=[x,y,p_1,\dots, p_N,z_1 \dots z_N, u, v ]^t$ as follows:
\begin{equation}\label{clari2}
\begin{array}{cccccc}
\dot x & = & 0 & +y  &  \\
\dot y & =& 0  & +\omega_r^2 p - q_r \omega_r y -\omega_r^2 x & \\
\dot p_n & = & 0 & +z_n & & \quad \forall n \in [1 \dots N] \\
\dot z_n - \frac{2c}{l} \dot u & = & 0 & -2\alpha_n c \, z_n -\omega_n^2 p_n &  & \quad \forall n \in [1 \dots N] \\ 
0 & = & 0 & -u + \zeta (1-\lambda) v & + \zeta x v \\ 
0 & = & 0 & -p + p_1 + \dots +p_N &  \\
 \underbrace{\qquad 0 \qquad }_{m(\dot Z)} & =& \underbrace{\quad \lambda \quad}_{c(\lambda)}  & \underbrace{\quad \quad \quad \quad -p \quad \quad \quad \quad}_{l(Z,\lambda)} &\underbrace{\quad -v^2   \quad }_{q(Z,Z)}   \\
\end{array}
\end{equation}
Here, the linear term $l$, but also the constant term $c$,  are functions of $\lambda$.

\subsection{The harmonic balance method applied to a quadratic system}

In this section, the harmonic balance method is applied to the system (\ref{equaquad}).
The unknown (column) vector $Z$ is decomposed into  Fourier series with $H$ harmonics:
\begin{equation}\label{fourierserZ}
Z(t)= Z_0 + \sum_{k=1}^H Z_{c,k} \cos(k \omega t) + \sum_{k=1}^H Z_{s,k} \sin(k \omega t)
\end{equation}
The components of the Fourier series  are collected into a large (column) vector $U$, with size $(2H+1)\times N_e$, where  $N_e$ is the number of equations in (\ref{equaquad}).
\begin{equation}\label{vectU}
U=  [ Z_0^t ,  Z_{c,1}^t , Z_{s,1}^t , Z_{c,2}^t , Z_{s,2}^t , \dots , Z_{c,H}^t , Z_{s,H}^t ]^t
\end{equation}
By introducing the expansion (\ref{fourierserZ}) into the set of Eqs. (\ref{equaquad}),  collecting the terms of the same harmonic index, and neglecting the higher order harmonics, one obtains a large system of $(2H+1) \times N_e$ equations for the unknown vector $U$, 
\begin{equation}\label{grandsys}
\omega M(U)= C + L(U) + Q(U,U)
\end{equation}
The new operators $M(.)$, $C$, $L(.)$, and $Q(.,.)$ that apply to $U$  depend only on the operators $m(.)$, $c$,  $l(.)$ and $q(.,.)$ of Eq. (\ref{equaquad}) and on the number of harmonics H.  The explicit formulas  have been given in annex 1 for a sake of brevity. It should be noted that these expressions for $M(.)$, $L(.)$ and $Q(.,.)$ are the same in all the examples presented. The change from one particular system to another only involves the operators  $m(.)$, $l(.)$ and $q(.,.)$.

The final system (\ref{grandsys}) contains $(2H+1) \times N_e$ equations for the $(2H+1) \times N_e$ unknowns $U$ plus the angular frequency $\omega$ and the continuation parameter $\lambda$.
Since the original system is autonomous, a phase condition has to be added to Eq. (\ref{grandsys}) to define a unique orbit. Indeed, the  time $t$ does not appear in Eq. (\ref{equaquad}) and if $U(t)$ is a solution of Eq. (\ref{grandsys}) then $U(t+\tau)$ is also a solution, for any $\tau$. We refer here to textbooks dealing with periodic solution continuation (\cite{Seydel94,Doedel07,nayfeh95}) for the choice of the phase condition. 

\subsection{The continuation procedure}

\subsubsection{Framework}

As the result of the harmonic balance operation, we  now have to solve an algebraic system 
\begin{equation}
\br(\bu)=0
\label{e:ru}
\end{equation}
 with  $\br \in \mathbb{R}^{N}$ and $\bu = [ U^t , \lambda, \omega ]^t \in \mathbb{R}^{N+1}$ ($N=(2H+1) \times N_e+1$). Any numerical continuation method (\cite{allgower90,krauskopf07,Seydel94}) and any suitable software (see (\cite{govaerts07} for an overview) can be used for this purpose. However, in this paper we will use the so-called Asymptotic-Numerical-Method (ANM), on which the idea of the quadratic recasting was based. It is worth noting that the final  system (\ref{grandsys}) is already quadratic with respect to  $U$ and $\omega$,
and that the application of the ANM is therefore quite straighforward. 

In the continuation procedure,  a pseudo-arc length parametrization will be used in order to be able to pass the limit points with respect to $\lambda$. The parameter $\lambda$ will therefore now become an unknown, like $U$ and $\omega$.  It is now mandatory to   specify how the operators $m$, $c$, $l$ and $q$ in Eq. (\ref{equaquad}) depend on $\lambda$. 
In the three examples presented here, we have organised the terms of equations so that the parmeter $\lambda$ only appears in the operators $c$ and $l$. In addition, we have managed to make this dependence linear. Once again, this formulation is not limited to the three example presented here, but it can be obtained for a very large class of dynamical systems  provided  suitable additional variables are introduced. We recall that the parameter $\lambda$ should not appear in the operator $q$, for which the HBM algebra is the more complex. For example, if there is a term  $\lambda u^2$ in  Eq. (\ref{equaquad}), it should be rewritten as $\lambda v$ with the auxiliary variable $v=u^2$, and put into $l$ instead of $q$.

In what follows, it is therefore  assumed that $c$ and $l$ can be written:

\begin{equation}\label{lambda}
\begin{array}{ll}
c & = c_0 + \lambda c_1 \\
l(.) &= l_0(.) + \lambda l_1(.)
\end{array}
\end{equation}
where $c_0$, $c_1$, $l_0$ and $l_1$ are independent of $\lambda$.
 Under this assumption, the operators $C$ and $L$ simply become
\begin{equation}\label{lambda2}
\begin{array}{ll}
C & = C_0 + \lambda C_1 \\
L(.) &= L_0(.) + \lambda L_1(.)
\end{array}
\end{equation}
Note that the expression for $C_0$, $C_1$ and $L_0$, $L_1$ in terms of $c_0$, $c_1$ and $l_0$, $l_1$ are exactely the same as those given in annexe 1 for  $C$ and $L$ in terms of $c$ and $l$.
The final algebraic system (\ref{e:ru}) becomes
\begin{equation}\label{sysman}
\br(\bu)=\bl0 \; + \; \bl(\bu) \; + \; \bq(\bu,\bu)
\end{equation}
with $\bu=[ U^t , \lambda, \omega ]^t $ and
\begin{equation}\label{sysman2}
\begin{array}{ll}
\bl0 & = C_0  \\
\bl(\bu) & = L_0(U) + \lambda C_1   \\
\bq(\bu,\bu) & = Q(U,U) \, + \, \lambda L_1(U)\,  - \, \omega M(U)  \\
\end{array}
\end{equation}
In Eq. (\ref{sysman}), $\bl0$ is a constant vector, $\bl(.)$ is a linear vector valued operator and $\bq(.,.)$ is a bilinear vector valued operator. We will now  briefly review  the ANM continuation technique for solving  quadratic system Eq(\ref{sysman}).

\subsubsection{The ANM continuation}

One of the main particularities of the ANM is that it gives access to branches of solution
in the form of power series. Assuming that we  know a regular solution point  $\bu_0$, the  branch of solution passing through this point is computed in the form of
a power series expansion (truncated at order $n$) of the pseudo-arclength path parameter $a=(\bu-\bu_0)^t \bu_1$, where  $\bu_1$ is the tangent vector at $\bu_0$:
\begin{equation}\label{U_serie}
\bu(a) = \bu_0 + a \bu_1 + a^2 \bu_2 + a^3 \bu_3 + \dots + a^n \bu_n.
\end{equation}
The series (\ref{U_serie}) is replaced in Eq. (\ref{sysman}) and each power of $a$ is equated to zero, giving a series of linear systems :

\begin{itemize}
\item {\bf order 0 :} $\bl0 + \bl(\bu_0) + \bq(\bu_0,\bu_0) = 0 $, which is obvious since $\bu_0$ is a solution of  Eq.  (\ref{sysman}).
\item {\bf order 1 :} $\bl(\bu_1) + \bq(\bu_0,\bu_1) +  \bq(\bu_1,\bu_0) = 0 $, which can also be written $\bj_{\bu_0}\bu_1 = 0$ where $\bj_{\bu_0} \in \mathbb{R}^{N \times N+1}$ is the jacobian matrix of $\br$ evaluated at $\bu_0$.
\item {\bf order $\bf 2 \leq p \leq n$ :} $\bj_{\bu_0}\bu_p + \Sigma_{i=1}^{p-1}\bq(\bu_i,\bu_{p-i})= 0$
\end{itemize}
The original nonlinear problem has therefore been reduced to a series of $n$ linear systems of $N_e$ equations. However, at each order, the linear systems are under-dimensioned since they have $N_e+1$ unknowns. The additional equations required are obtained by inserting Eq. (\ref{U_serie}) into the definition of the path parameter $a$ given above.
This gives:
\begin{itemize}
\item {\bf order 1 :} $\bu_1^t  \bu_1 = 1$
\item {\bf order $\bf 2 \leq p \leq n$ :} $\bu_1^t \bu_p = 0$
\end{itemize}

\subsubsection{Comments:} 
\begin{itemize}
\item The original nonlinear system of $N_e$ equations has been
  transformed into $n$ linear systems of $N_e$ equations, which have to be solved
  successively (as in classical perturbation methods), i.e., $\bu_p$ is deduced from terms occuring at lower orders.
\item The range of utility of a truncated power series is generally limited because the series have a finite radius of convergence.   Once each $\bu_p$ ($p \in [1 \dots N]$) has been found, the range of utility of the
  series expansion  (\ref{U_serie}) is defined by the value $a_{max}$ such that 
  \begin{equation}
  \forall a \in [0  \; \am], ||\br(\bu(a))||  \leq \epsilon_r 
\label{e:amax} 
 \end{equation}
where $\epsilon_r$ is a user-defined tolerance parameter (see (\cite{cochelin94c,manlivre}) for details of the $a_{max}$ calculation). Note that the range of utility $a_{max}$  is generally approximately equal  to the radius of convergence of the series (\cite{baguet03}).
\item The series expansion and the associated range of utility $a_{max}$ only define one part
of the branch of solutions, which is called a section. To determine the entire branch, it is necessary to 
restart the whole series calculation from successively updated starting points $\bu_0$. The simplest way of performing this updating is to take the end point of a previously computed section as the starting point for the next section. This is the principle of continuation with the ANM.
Finally, the entire branch is given
  as a succession of different  branch sections, where the length of each
 section is automatically given by its range of utility $a_{max}$. 
\item The ANM continuation approach has the following advantages:
\begin{itemize}
\item the solution branch is known analytically, section by section.
\item since all the linear systems to be  solved have the same Jacobian matrix, the computational cost of the series (\ref{U_serie}) is low.
\item since the size of the section is given by the convergence properties of the current step, ie, by the values of $a_{max}$, the ANM continuation algorithm does not require any special  step-length control strategies. 
\item the ANM continuation algorithm is highly robust, even when the branch contains sharp turns.  Branch switching can therefore be easily performed, using small  pertubations in the system.
\end{itemize}
\end{itemize}

\subsection{Implementation in the MANLAB software program}

\subsubsection{A brief presentation of MANLAB}
MANLAB is an interactive software program for the continuation and bifurcation analysis of 
algebraic systems, based on ANM continuation. The latest version is programmed in Matlab using an objet-oriented approach
(\cite{remithese}). 
MANLAB has a Graphical User Interface (GUI) with buttons, on-line inputs and graphical windows  for generating, displaying and analysing the bifurcation diagram and the solution of the system. To enter the system of equations, the user has to provide three vector valued Matlab functions corresponding to the constant, linear and quadratic operators $\bl0$, $\bl$ and $\bq$.
As an  example, take the biochemical reaction system used by Doedel et Al in (\cite{Doedel91})
 \begin{equation}
\label{eq:exemple}
\begin{array} {lll}
r_1(u_1,u_2,\lambda)= & 2 u_1 -u_2 + 100 \frac{u_1}{1+u_1+u_1^2}-\lambda &=0 \\
r_2(u_1,u_2,\lambda)= & 2 u_2 -u_1 + 100 \frac{u_2}{1+u_2+u_2^2}-(\lambda+\mu)& =0.
\end{array}
\end{equation}
Introducing the following additional variables $v_1  = u_1 + u_1^2$,  $v_2  = u_2 + u_2^2$,  $v_3  = \frac{1}{1+v_1}$
and  $v_4  = \frac{1}{1+v_2}$, the system can be rewritten with  quadratic polynomial nonlinearities as follows,
\begin{equation}
\label{eq:exqua}
\begin{array}{lclclll}
  0    &  + \quad & 2u_1-u_2 -\lambda  & \quad + \quad & 100 u_1 v_3 &=  0 \\
 -\mu &  + \quad & 2u_2-u_1-\lambda   & \quad + \quad & 100 u_2 v_4  &=  0 \\
  0    &  + \quad & v_1 - u_1          & \quad -\quad  & u_1^2       & =  0 \\
  0    &  + \quad  &v_2 - u_2          & \quad - \quad & u_2^2       & =  0 \\
 -1   &    +\quad & v_3               &\quad  + \quad & v_1 v_3     &=  0 \\
\underbrace{ -1  }_{\bl0} &  + \quad   &  \underbrace{ v_4 \quad \quad \quad \quad  \quad}_{\bl(\bu)}   &\quad + \quad& \underbrace{\quad v_2 v_4 \quad}_{\bq(\bu,\bu)}    & =  0
\end{array}
\end{equation}
The first two equations  correspond directly to Eq ( \ref{eq:exemple}). The last four define the auxilliary variables $v_1, v_2, v_3, v_4$. The constant, linear and quadratic terms in each equation have been split to  define the three operators $\bl0$, $\bl$ and $\bq$. 
In the MANLAB software program, the unknowns are collected in a single vector $\bu =[ u_1, u_2, v_1, v_2, v_3, v_4, \lambda]$, and the three operators are given by the following vector valued functions
( $\mu=0.05$ in the case of this example ) :

{\scriptsize
\begin{verbatim}
function [L0] = L0           function [L] = L(U)                 function [Q] = Q(U,V) 
                      
L0=zeros(6,1);               L=zeros(6,1);                       Q=zeros(6,1); 
L0(1)=   0;                  L(1)=2*U(1)- U(2) -U(7);            Q(1)=100*U(1)*V(5); 
L0(2)=  -0.05;               L(2)=2*U(2)- U(1) -U(7);            Q(2)=100*U(2)*V(6);
L0(3)=   0;                  L(3)=U(3)-U(1);                     Q(3)=   -U(1)*V(1);
L0(4)=   0;                  L(4)=U(4)-U(2);                     Q(4)=   -U(2)*V(2);
L0(5)=  -1;                  L(5)=U(5);                          Q(5)=    U(3)*V(5);
L0(6)=  -1;                  L(6)=U(6);                          Q(6)=    U(4)*V(6);
\end{verbatim}
}

\subsubsection{Implementation of the periodic solution continuation in MANLAB}

To compute the branches of periodic solutions, the system has to be first recasted  in the form of Eq. (\ref{equaquad}) with account of the additional decomposition (\ref{lambda}). This is probably the most unusual and difficult task for a new user. Subsequently, it is only necessary to provide the MANLAB software program with 
the operator $m(.)$, $c_0$, $c_1$, $l(.)$ $l_1(.)$ and $q(.,.)$ . The functions  $\bl0$, $\bl$ and $\bq$, which are the actual input for MANLAB, have been programmed once for all, using the expression  in Eq. (\ref{sysman2}) and the formulas given in annex 1. The examples presented in this paper are available online (\cite{SiteManlab}). 

\section{The case of a periodically forced system \label{s:forced}}

We now focus on periodically forced (non-autonomous) systems:
\begin{equation}\label{systdynana}
\dot Y=f(t,Y,\lambda)
\end{equation}
where $f$ is periodic in $t$, with the period $T$ (forcing period). We look for periodic solutions (responses) with
a period $pT$ or $\frac{T}{p}$, where $p$ is an integer.  A classical strategy consists in  transforming Eq. (\ref{systdynana}) into an augmented autonomous system (\cite{Seydel94}),   by adding an oscillator with the desired forcing period in the system of equations, for instance (\cite{Doedel07}).

In what follows, a direct approach will be used. It consists in expanding the forcing term into harmonics, and taking them into account in the balance of individual harmonics. This appraoch is illustrated below with two further examples.

\medskip

\textbf{Example 4: } Forced Duffing oscillator.

The normalized forced Duffing oscillator is the non-autonomous equation:
\begin{equation}\label{duf}
\ddot u +2  \mu \dot u + u  + u^3 = f \cos( \lambda t)
\end{equation}
We take the damping coefficient  $\mu$ and the force amplitude $f$ constant, and  use the forcing angular frequency as the varying parameter $\lambda$.
By using $v(t)= \dot u(t)$ and $w(t)= u^2(t)$, this equation can be recast as follows
\begin{equation}\label{duf2}
\begin{array}{ccccc}
\dot u & = & & v &  \\
\dot v & = & f  \cos(\lambda t)  & -2 \mu v - u   &  - u w \\
 \underbrace{0 }_{m(\dot Z)} & = & \underbrace{\quad \quad 0\quad \quad}_{c(t, \lambda)}   & \underbrace{\quad   \quad w \quad \quad   }_{l(Z)}  \  & \underbrace{\quad  \quad  - u u  \quad \quad }_{q(Z,Z)}  \\
\end{array}
\end{equation}
where $Z=[ u , v , w]^t$, and the forcing term is deliberately put into $c$.
The forcing frequency is now related to  the response frequency by putting $\lambda=\omega$ or possibly, $\lambda=p \omega$ (p is an integer). The term $c(t)$ is then expanded into harmonics with respect to $\omega$.

This results in  slight changes  in the procedure:
\begin{itemize}
 \item because of the synchronization of the response and the forcing, the phase condition has to be removed. Note that the parameter $\lambda$ is no longer an unknown, since it was chosen as a multiple of $\omega$. In comparison with  the case of an autonomous system, both the number of equations and the number of unknowns have decreased by one.
 \item in the final system (\ref{sysman})(\ref{sysman2}), the operators $C_1$ and $L_1$ disappear and  the forcing amplitude $f$ has to be accounted for in $L_0$ 
\end{itemize}

Lastly, we take the Raylegh-Plesset equation, which is used to model the large amplitude vibrations of a gas bubble in a fluid. The  forcing is handled slightly differently and this example also shows  how to cope with a power of $-3$.

\medskip 

\textbf{Example 5: } The  forced Rayleigh-Plesset equation (\cite{Plesset_1949}). 

Let $R$ be the radius of the vibrating bubble, and $R_0$ the radius at rest. The equation of motion of the bubble is
\begin{equation}\label{RP}
R \ddot R + \frac{3}{2} \dot R^2 =  A  \{ (\frac{R_0}{R})^3 -1 \} + B \cos( \lambda t)
\end{equation}
where $A$ and $B$  are fixed. 
We introduce the normalized radius $u=\frac{R}{R_0}$ and the (normalized) velocity $v= \frac{\dot R}{R_0}$.
Dividing Eq (\ref{RP}) by $R_0^2$ and defining $a=\frac{A}{R_0^2}$, $b=\frac{B}{R_0^2}$, we get the first order system
\begin{equation}\label{RP2}
\begin{array}{lll}
 \dot u & =  & v \\
 u \dot v  & =  &   a (u^{-3}-1)   -\frac{3}{2} v^2 + b \cos(\lambda t)  \\
\end{array}
\end{equation}
We now introduce the following auxiliary variables $x=\frac{1}{u}$,  $y=x^2$, and $z=v^2$, and arrive at 
\begin{equation}\label{RP3}
\begin{array}{lll}
 \dot u & =  & v \\
 \dot v  & =  &   a (y^2 - x)   -\frac{3}{2} x z + b x \cos(\lambda t)  \\
\end{array}
\end{equation}
After undergoing this transformation, the system is quadratic but  the forcing term has now been  multiplied by the unknown function $x$, and it is no longer a constant term with respect to the unknown. We introduce another auxilliary variable $r(t)=  \cos(\lambda t)$, and replace the term $b x \cos(\lambda t)$ by $b x r$. Finally, we obtain the following system, where the first two equations stand for Eq (\ref{RP3}) and the  last four define the auxilliary variables.
\begin{equation}\label{RPsys}
\begin{array}{ccccccc}
\dot u & =  & 0 &\quad + \quad & v   &\quad + \quad  &  \\
\dot v & =  & 0 &\quad + \quad &(- a x) &\quad + \quad  &  a y^2 -\frac{3}{2} xz + b x r  \\
     0 &  = & 1 &\quad + \quad & 0   &\quad + \quad  & (- u x) \\
     0 & =  & 0 &\quad + \quad & y   &\quad + \quad  & (- x x) \\
     0 & =  & 0 &\quad + \quad & z   &\quad + \quad  & (- v v)  \\
\underbrace{0 }_{m(\dot Z)} & = & \underbrace{ -\cos(\lambda t) }_{c(t)}   &\quad + \quad & \underbrace{  \quad \quad r \quad \quad }_{l(Z)}   &\quad + \quad& \underbrace{    \quad  \quad   \quad \quad  \quad   \quad   \quad  \quad   \quad}_{q(Z,Z)}  \\
\end{array}
\end{equation}
Here the unknown vector is  $Z=[ u , v , x, y, z, r ]^t$, and the forcing term is clearly present in the operator $c$.

\section{Numerical results on selected examples}

In the  following selected examples, the numerical results obtained using the approach presented in this paper are either compared with time-domain simulations to show the validity of our approach, or used to illustrate particular features: the influence of the number of harmonics in the case of the Van der Pol oscillator, the ability to follow period-doubling bifurcations in that of the Rössler system, the ability to follow a direct or inverse Hopf bifurcation in that of the clarinet model, and to illustrate a forced system, in that of the Duffing oscillator.

\begin{figure}[!h]
\begin{center}
\begin{tabular}{cc}
\hspace*{-1.25cm}\includegraphics[width=0.6\textwidth]{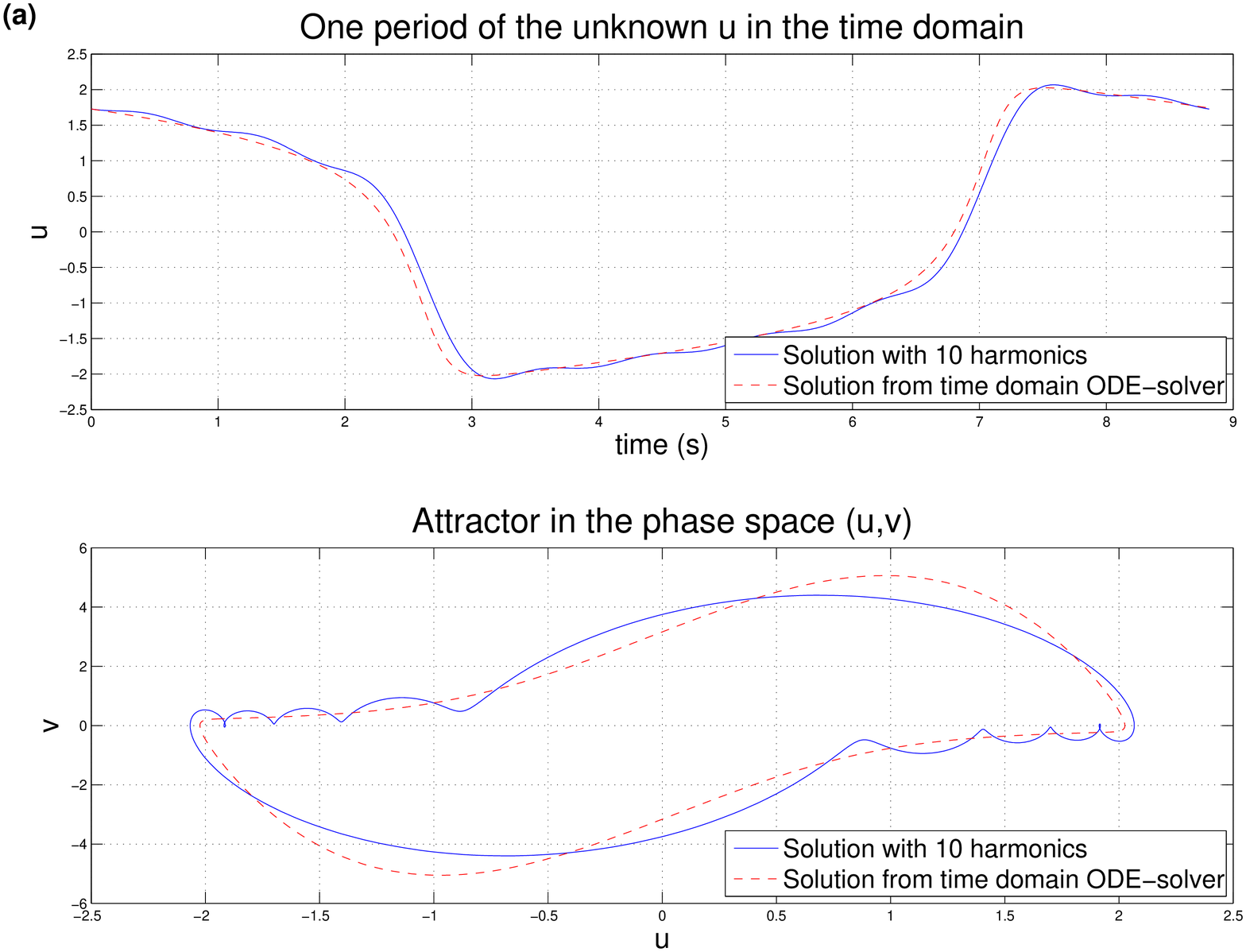}&
\hspace*{-1.25cm}\includegraphics[width=0.6\textwidth]{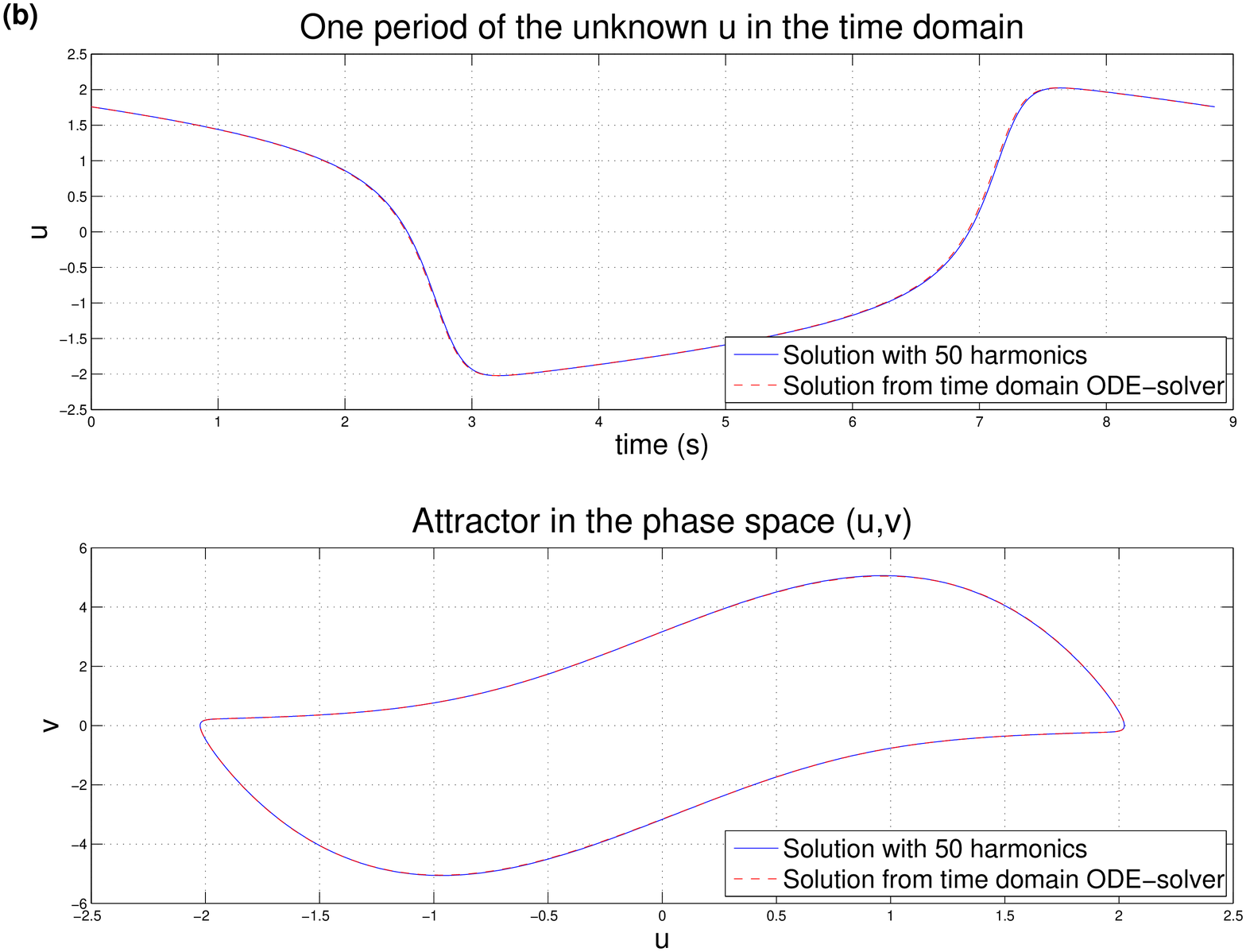}
\end{tabular}
\end{center}
\caption{The Van der Pol oscillator (Eq. (\ref{vdp1})): Comparison between  the harmonic balance results and those obtained  with a time domain simulation. Influence of the number of harmonics $H$ taken into account ((a) H=10, (b) H=50).\label{f:vdp}}
\end{figure}

\textbf{Example 1: The Van der Pol oscillator.}
Numerical results obtained on the Van der Pol oscillator (Eq. (\ref{vdp1})) are presented in Fig. \ref{f:vdp}.  Good agreement was obtained with the results of a time domain simulation, performed using Matlab ODE solvers. As expected, the number of harmonics $H$ in the solution sought by the harmonic balance method had to  be adapted to span the bandwith of the solution obtained by direct integration. In the case of Fig. \ref{f:vdp}, $\lambda=3$, which requires having a relatively high $H$ to match the reference solution.

\medskip

\begin{figure}
\includegraphics[width=0.75\textwidth]{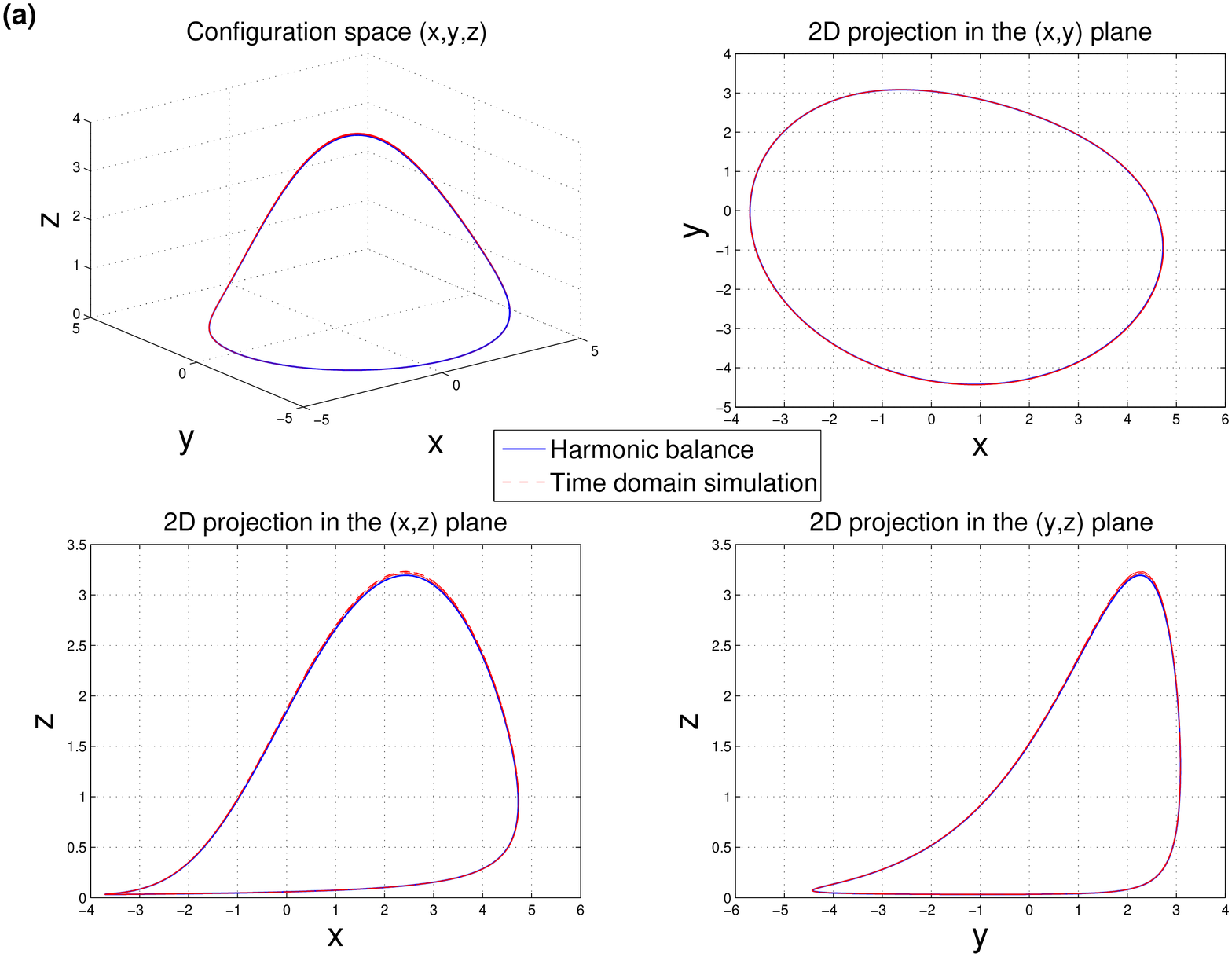} \vspace*{-0.25cm}
\includegraphics[width=0.75\textwidth]{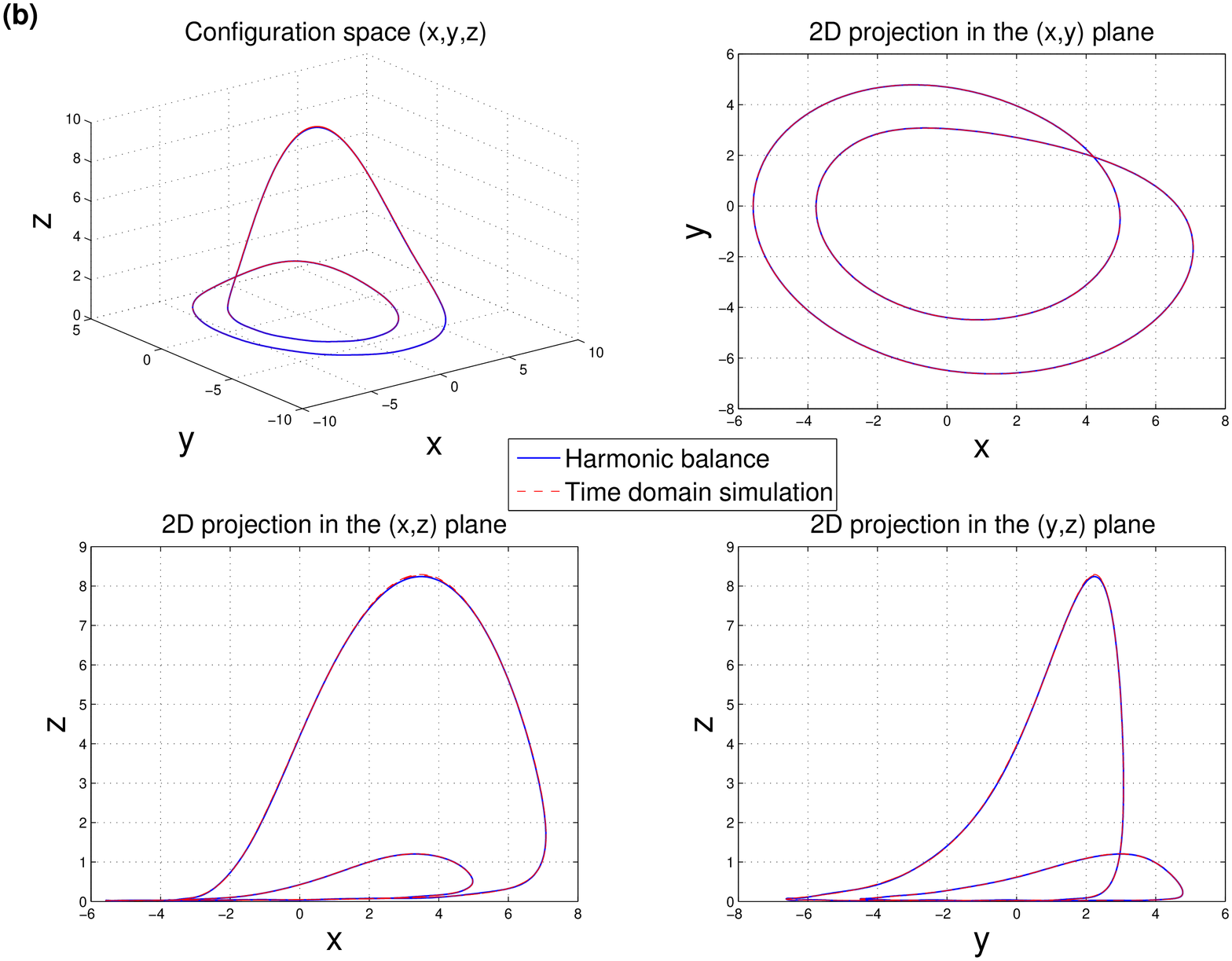} \vspace*{-0.25cm}
\includegraphics[width=0.75\textwidth]{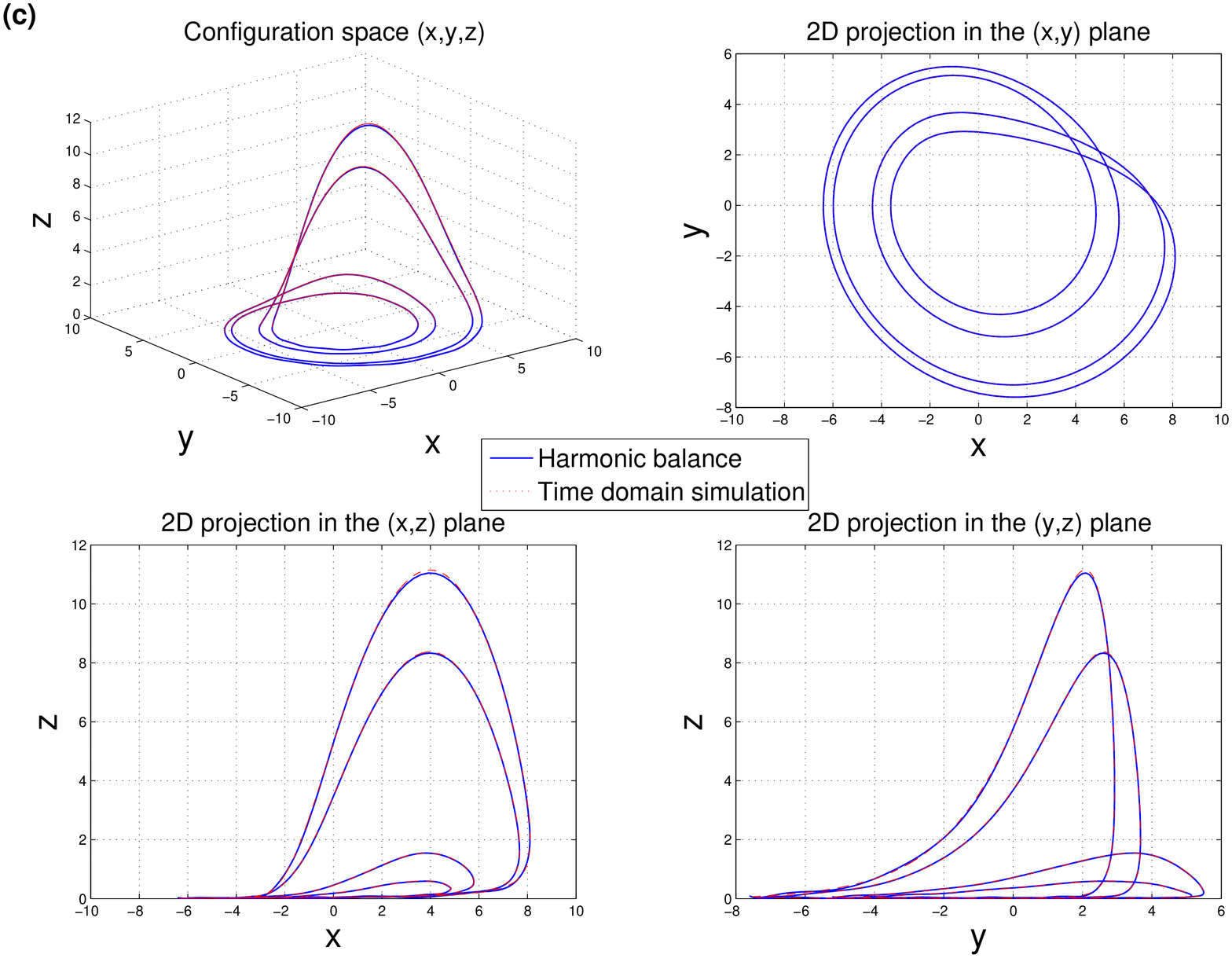} \vspace*{-0.25cm}
\caption{The Rössler system (Eq. (\ref{rossler1}) with $a=b=0.2$): Comparison between the results obtained using the harmonic balance method and a time domain simulation. These three figures show the ability of the method to follow the sub-harmonic cascade from (a) $\lambda=2.5$ (perdiod $T$), to  (b) $\lambda=3.5$ (period $2T$ after a first period-doubling bifurcation) and (c) $\lambda=4$ (period $4T$ after a second period-doubling bifurcation) \label{f:rossler}}
\end{figure}

\textbf{Example 2: The Rössler system.} 
The ability to follow period-doubling bifurcations is illustrated in Fig. \ref{f:rossler}. Following these bifurcations is more difficult in the case of autonomous systems than in forced systems, since the period of oscillation is also  unknown. Moreover, since a $T$-periodic solution is also a $2T$-periodic solution, one cannot expect the value of the period obtained  by  the harmonic balance process to be doubled when crossing a period-doubling bifurcation.  The strategy used here therefore consists in introducing  complex subharmonic amplitudes as additional unknowns. To detect $K$ period-doubling bifurcations, the  solution is then sought as:
\begin{equation}\label{fourierserZsubharm}
Z(t)= Z_0 + \sum_{k=1}^H Z_{c,k} \cos(\frac{k}{2^K} \omega t) + \sum_{k=1}^H Z_{s,k} \sin(\frac{k}{2^K} \omega t)
\end{equation}
When the solution belongs to the  $T$-periodic solution branch, $\left.Z_{c,k}\right|_{k=1..K-1}=0$ and $\left.Z_{s,k}\right|_{k=1..K-1}=0$.  When the solution belongs to  the $2T$-periodic solution branch, $\left.Z_{c,k}\right|_{k=1..K-1, k\neq2^{K-1}}=0$ and $\left.Z_{s,k}\right|_{k=1..K-1, k\neq2^{K-1}}=0$.  One practical consequence in terms of  the computational cost is that in order to span the same bandwidth, the number of harmonics has be mulitplied by $2^K$. This case is illustrated in Fig. \ref{f:rossler} where $H=10$ from the $T$-periodic solution to $H=2^2 \times 10 =40$ after two period-doubling bifurcations.

\medskip

\textbf{Example 3: The clarinet model.} 

\begin{figure}
\includegraphics[width=0.5\textwidth]{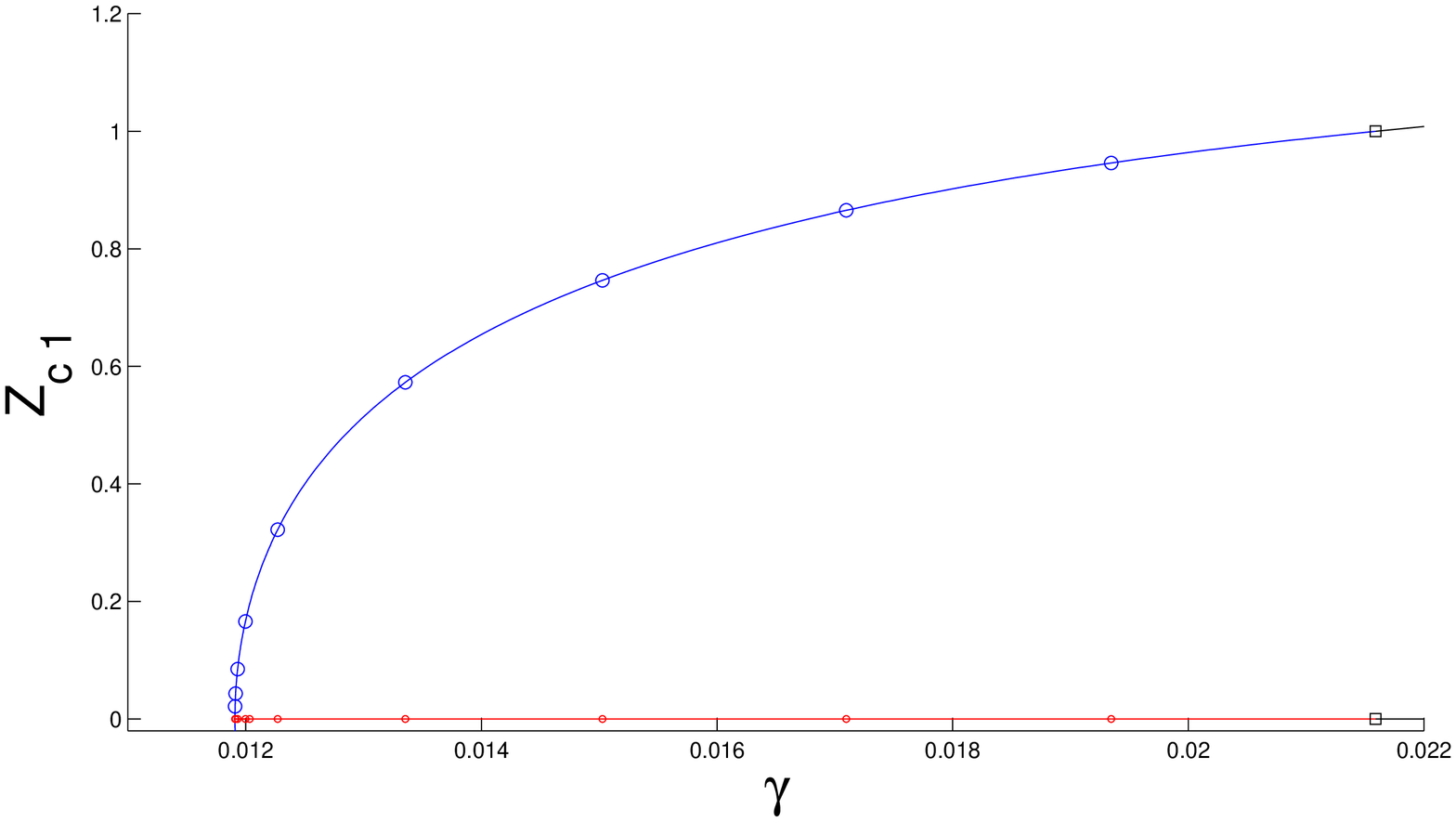} 
\includegraphics[width=0.5\textwidth]{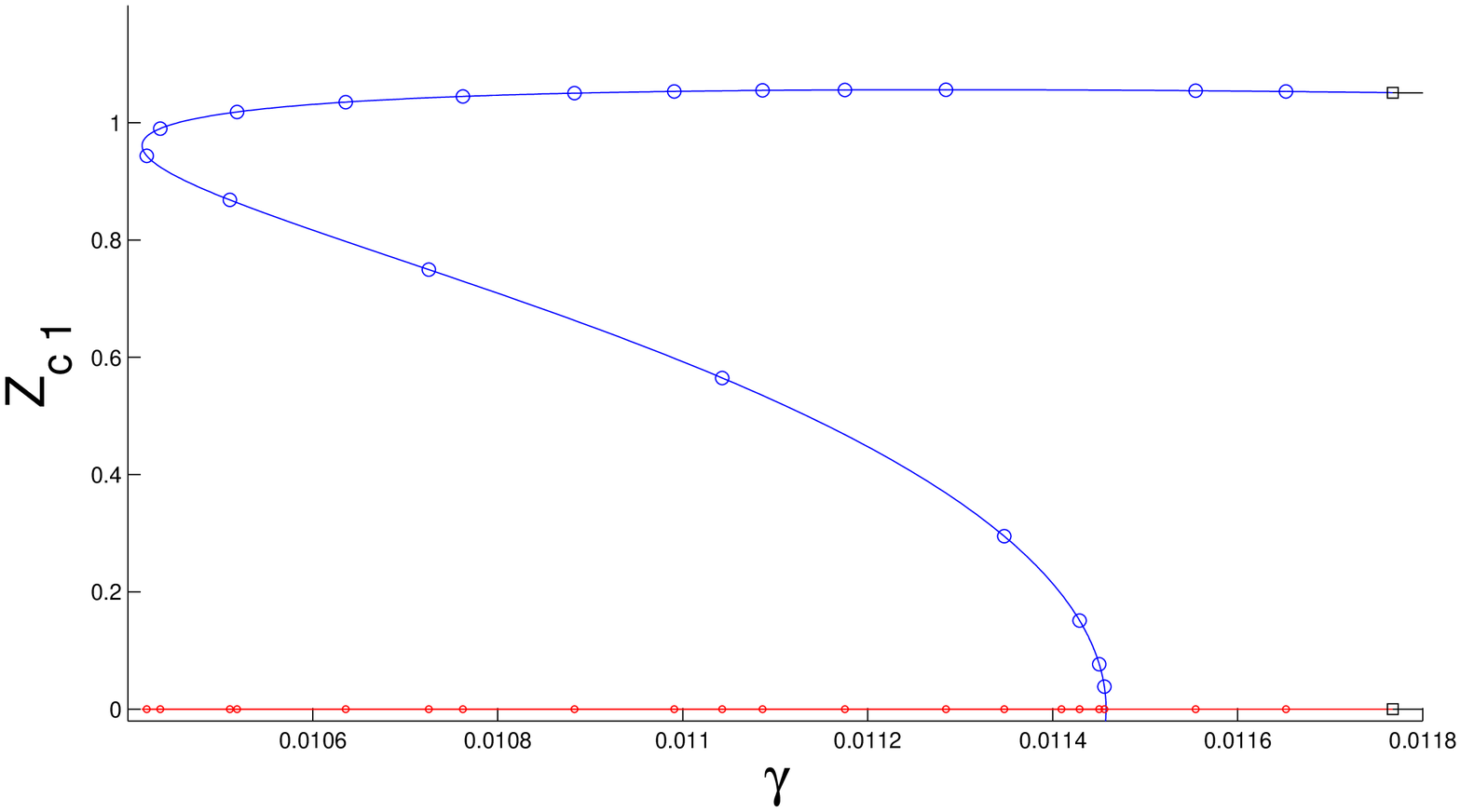}
\caption{The clarinet model (Eq. (\ref{clari1}) with $q_r=0.01$, $\omega_r=6597$, $H=3$, $N=5$, $c=340$,  $\zeta=0.2$, $l_t=5.6e-08$, $l_v=4e-08$, $C_{pv}=1.4$, $r=0.007$, $w_1=2426.5$, 
$w_2= 7279.5$, $w_3=12132.5$, $w_4=16985.4$, $w_5=21838.4$, $\alpha_1=0.11$, $\alpha_2=0.19$, $\alpha_3=0.25$, $\alpha_4=0.3$, $\alpha_5=0.34$): depending on the length of the resonator, the bifurcation will be either a direct Hopf bifurcation (left picture, $L=0.22$) or an inverse Hopf bifurcation (right picture, $L=0.2$).  \label{f:clarinet}}
\end{figure}

Two typical examples of the bifurcation diagrams obtained with the clarinet model are shown in Fig. \ref{f:clarinet}. These pictures focus on the oscillation threshold in order to show the robustness of the method, even at singular points. As can be seen from both pictures, the radius of convergence of the power series expansion decreases  when approaching a bifurcation point (and hence the length of each section of a branch lying between two consecutive points in Fig. \ref{f:clarinet} decreases). This behaviour is typical of ANM (\cite{baguet03}).

\medskip

\textbf{Example 4: The forced Duffing oscillator} 

Figure 4 shows the frequency-amplitude diagram of the response obtained with the  method presented here. A classical bent resonance curve can be observed, as well as some additional peaks corresponding to superharmonic resonances. Note that only the individual amplitude $A_i=\sqrt{(u_{c \; i}^2+u_{s \; i}^2)}$ of the odd harmonics have been plotted, since the even harmonics ($A_0$, $A_2$, $A_4$, $\dots$) are zero.  The computation of this branch of periodic orbits required 25 steps of MAN-continuation when 5 harmonics were included, and 35 steps when 9 harmonics were included. The curves obtained for  harmonics 1 and 3 were slightly hanged when shifting from $H=5$  to $H=9$ harmonics in  Eq. (\ref{fourierserZ}). The curve obtained for $A_5$ is more affected, which confirms the logical result that more than 5 harmonics have to be used to obtain an accurate result with  harmonic $5$. We  also checked that the curve obtained for $A_5$ was only very slightly modified when shifting from $H=9$  to $H=11$ harmonics in  Eq. (\ref{fourierserZ}).

\begin{figure}
\includegraphics[width=0.49\textwidth]{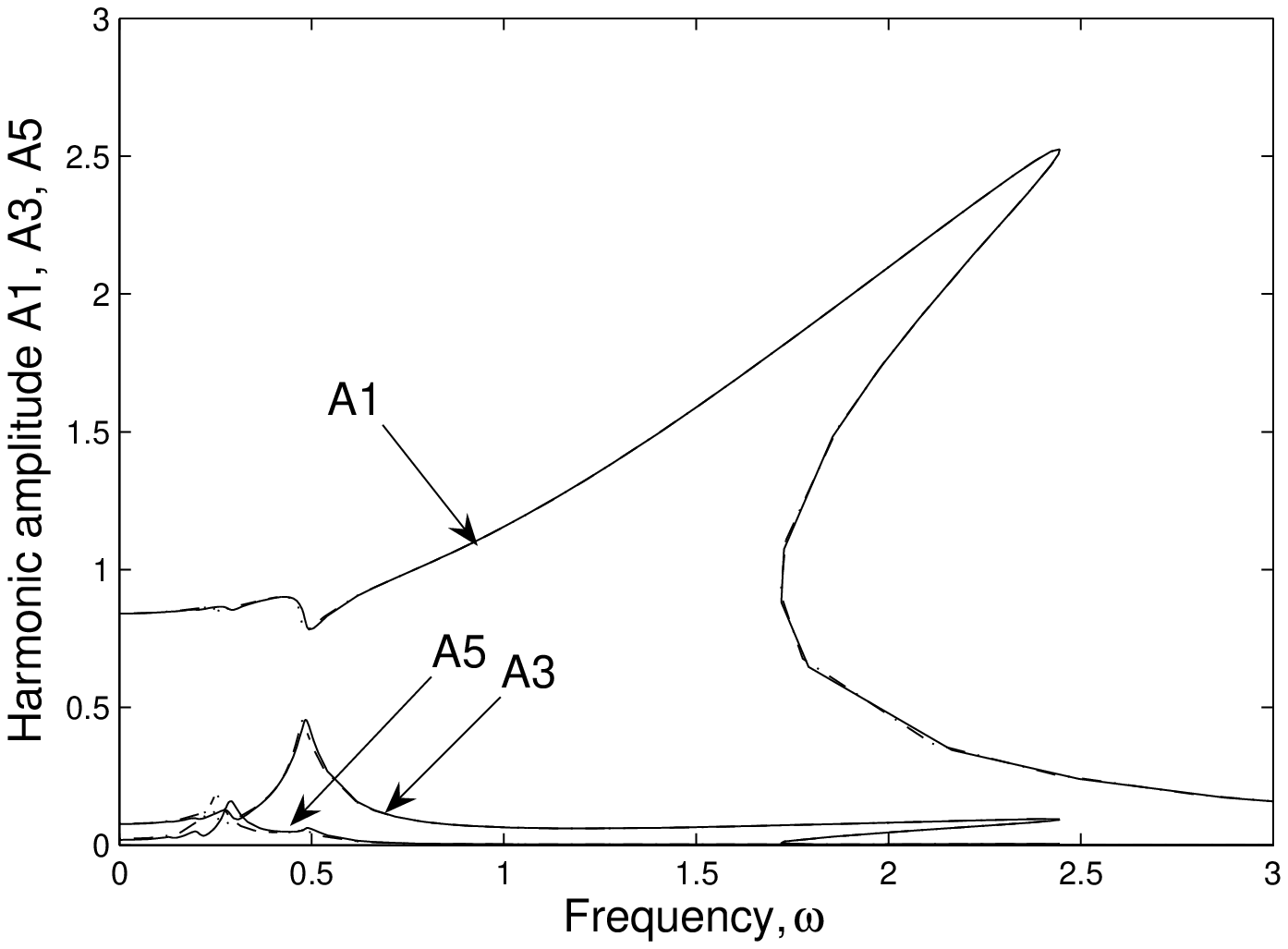}
\includegraphics[width=0.49\textwidth]{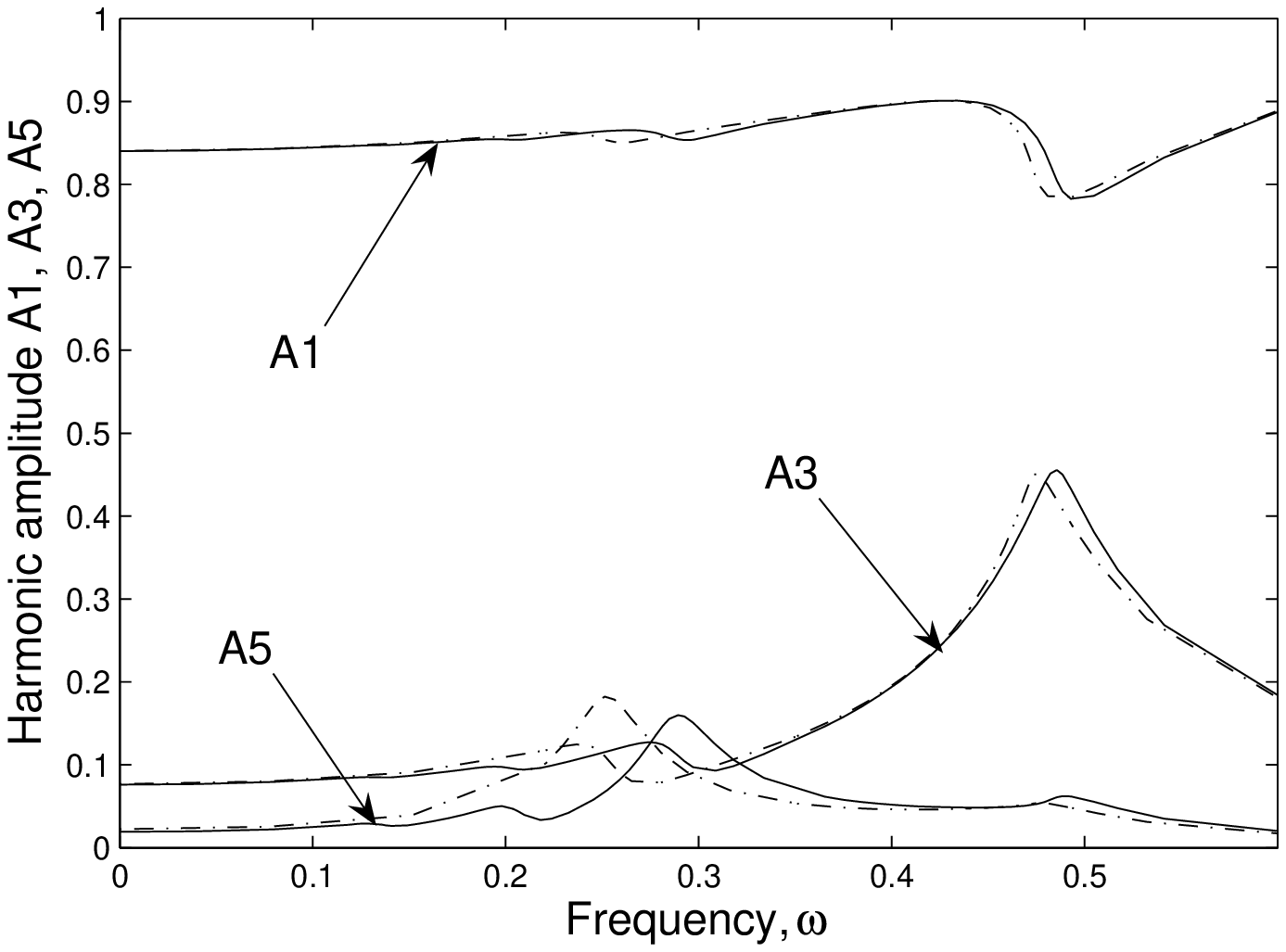}
\caption{Duffing oscillator, Eq. (\ref{duf}) with $\mu=0.1$ and $f=1.25$ as in (\cite{liu2006}). The figure shows the amplitude of the harmonics 1, 3 and 5 versus the frequency forcing $\omega$. Solid line : results when 5 harmonics are included (H=5 in Eq. (\ref{fourierserZ})). Dash line : results when 9 harmonics are included. The figure on the right is a zoom of the region $\omega < 0.6$. \label{f:duffing}
}
\end{figure}

\section{Discussion and Conclusion}

The key idea in this study is the quadratic recast. It is worth pointing out that the present method contains two quadratic recasts. The first one, which occurs in the time domain,  transforms the original system (\ref{systdyna}) into (\ref{equaquad}), and is presented in section \ref{ss:quadratic_recast}. The second one transforms the algebraic system (\ref{grandsys}) obtained with the HBM, into the final system (\ref{sysman}) which is solved by the ANM.  The operators in (\ref{grandsys}) and (\ref{sysman}) are linked  by (\ref{lambda2}) and (\ref{sysman2}). This point accounts for  many of the characteristics of  the method presented here: 
\begin{itemize}
\item The whole procedure takes place in the frequency domain using an analytical expression for the system  to be solved  (as with the classical harmonic balance) and allows to use an arbitrarily high number of harmonics  (as with the AFT procedures).
\item In the ANM continuation, no iterative process are used to solve the algebraic system, and  hence no problems of convergence arise. The computational cost is almost predictible and fairly rather low. On the contrary, convergence of the iterative process can be difficult with AFT procedures (\cite[sect. 4.6, p117]{Troyanovsky_1997}).
\item There is no need for FFT/IFFT procedures, which increases the computational load in the case of AFT approaches.
\item No  temporal discretization has to be performed (since  the whole procedure is in the frequency domain), and there is therefore no need to cope with aliasing.
\item There is no need to use numerical schemes to compute derivatives. Since the problem is a quadratic one , the Jacobian is calculated analytically and this can be done very easily. This is a considerable advantage of the approach presented here. Indeed, as stated by \cite[p17]{Troyanovsky_1997}, usually "analytical calculation of the Jacobian involves considerable calculation and transformation" but on the other hand, when considering numerical calculation of the Jacobian, "for very nonlinear circuits, the error introduced in the Jacobian estimation results in inaccurate updates of the unknowns and a large number of iterations and possibly no convergence".
\end{itemize}

As stated in \cite{peng2008}, to apply HB, "it is always necessary to write specific computation programs for different nonlinear models". Here, on the contrary,  the procedure is simple and rather general because of the natural splitting between the automated part, which is independent of the problem (i.e. the transformation from operators $c, l, q$ to operators $\bl0, \bl, \bq$ and the ensuing solving using the ANM), and the part which depends on the particular problem to be solved (the writing of operators $c, l, q$) which is let to the users. From the point of view of the users, this means that the work required is relatively easy and the software is simple to use. In addition, very few parameters have to be defined by the users: the number of harmonics in the solution and  parameter $\epsilon_r$ in Eq. (\ref{e:amax}), which controls the accury of the ANM computation.

In this paper, we have presented only small-sized systems for the sake of simplicity. The method can be obviously extended in principle to large systems of equations  although this may require some additional technical work, and introduce limitations due to the lack of computer ressources.  In \cite{perignon02}, the method  has been applied to compute the forced responses of geometrically nonlinear elastic structures discretized using the finite element method.
It is well known that with problems of this kind, the equations of motion are second order differential equations in terms of the displacement nodal variables $u$, and that they show a polynomial cubic nonlinarity.  Choosing the nodal velocity $v$ and the second Piola-Kirchhoff $S$ (at the Gauss point) as  auxiliary variables, the system of equations can easily be written in the form of a first order system with quadratic polynomial non linearities.When a cosine forcing term is applied, they read (in the discrete form)
\begin{equation}\label{elasticmodel}
\begin{array}{lll}
\dot u & = & v \\
\mathbf M \dot v & = &-\int_V (B_L + B_{NL}(u))^t S dv  + \mathbf F cos(\lambda t) \\
S & = & D (B_L + \frac{1}{2} B_{NL}(u)) u
\end{array}
\end{equation}
where $\mathbf M$ is the mass matrix, $D$ is the Hooke operator and $B_L$ and  $B_{NL}$ are the linear and nonlinear strain-displacement operators \cite{crisfield91}. In this form, the governing equations consist in a group of first order differential equations (velocity definition, equation of motion) and a group of algebraic equations (constitutive law at the Gauss point of the finite elements). Both groups are quadratic in $u$ and $S$. The algebra given in annex 1 has been directely introduced into a home-made FEM code, at the element level.


To conclude, it is proposed to  discuss the limitations of this method and suggest some possible lines  of  future research.
First of all, it is worth noting that when the original system is not directly in the  same quadratic form as Eq. (\ref{equaquad}),  introducting auxiliary variables,  as done in this paper, examples may lead to additional (possibly non-physical) solutions. For example, in the case of the clarinet model (example 3), a variable $v=\sqrt{\gamma - p}$ is introduced, and then further defined by the quadratic equation $v^2 = \gamma - p$. This definition corresponds to $v=\pm \sqrt{\gamma - p}$. In practical terms, this means that the user has to check that the branch of solution computed corresponds to a positive $v$.

 The main limitation of the  method presented is that the original systems cannot all be put into the quadratic form  Eq. (\ref{equaquad}). For this reason, problems where one or several relations are given in the frequency domain (typically, when some variables are linked through an impedance) cannot be treated. Similar difficulties arise when the original system contains functions (trigonometric functions, exponential functions, non-integer powers function, etc) that are applied to the unknowns. The classical pendulum model $\ddot \theta + \sin(\theta) =0$ is a good example of  systems of this kind. In these cases, we can generally introduce new variables and new differential equations, so that the functions are generated by the system itself. In the case of  pendulum model,  $v=\dot \theta$ and the two new variables $s(t)=\sin(\theta(t))$ and $c(t)=\cos(\theta(t))$ are introduced. Taking the derivative with respect to time for the last two equations, we obtain
\begin{equation}
 \begin{array}{lll}
 \dot \theta & = & v \\
\dot  v & = & -s \\
\dot s & = & c v \\
\dot c & = & -s v
\end{array}
\label{pendulquad}
\end{equation}
This quadratic system corresponds to the pendulum model, provided  the following two initial conditions are added
\begin{equation}
 \begin{array}{lll}
 s(0) & = & \sin(\theta(0)) \\
c(0) & = & \cos(\theta(0)) 
\end{array}
\label{pendulinit}
\end{equation}
The application of the HBM to Eq. (\ref{pendulquad}) will now be quite straighforward, but, because of Eq. (\ref{pendulinit}), the final
algebraic system will not be quadratic. In these cases, the application of the ANM continuation is still feasible but this requires  more elaborate algebra for computating  the power series, as explained in  (\cite{potier97,manlivre}). 
The efficiency  of the present method on the pendulum model will be tested in a near future. Another interesting idea would be to test whether this method  can be used to solve regularized non-smooth dynamic problems, and for instance, the case of vibrating systems with contact conditions and friction laws.

\section*{Acknowledgements}
The authors want to warmly thank Marie-Christine Pauzin and Olivier Thomas for their suggestions and for testing the programs avaible at \verb!http://manlab.lma.cmrs-mrs.fr!.
The authors are grateful to Benjamin Ricaud for useful discussions about the present paper.

The research was supported by  French National Research Agency \textsc{anr} in the context of the \textsc{Consonnes} project.


\appendix

\section{Appendix}
Here, we give  the expression for the operators $M$, $C$, $L$ and $Q$ in Eq. (\ref{grandsys})
in terms of $m$, $c$, $l$, $q$ and $H$. We recall that Eq. (\ref{grandsys}) is obtained
by substituting Eq. (\ref{fourierserZ}) into Eq. (\ref{equaquad}) and by collecting the terms
with the same harmonic index (cosines and sines).
The inputs of $M$, $L$ and $Q$ are the vector $U$, which contains the Fourier coefficients of $Z(t)$
\begin{equation}\label{vectU2}
U=  [ Z_0^t ,  Z_{c,1}^t , Z_{s,1}^t , Z_{c,2}^t , Z_{s,2}^t , \dots , Z_{c,H}^t , Z_{s,H}^t ]^t
\end{equation}

\subsection{Constant term}

The constant vector $c$ is associated with the harmonic zero. The operator $C$ is simply
\begin{equation}\label{opeC}
C=  [ c^t ,  0^t , 0^t ,  \dots  ]^t
\end{equation}

\subsection{Linear term}

For the linear operator $l$, we have 
{\footnotesize
\begin{equation}\label{lin1}
\hspace*{-0.2cm}
\begin{array}{ll}
l(Z(t)) & = l(Z_0+\sum_{k=1}^H Z_{c,k} \cos(k \omega t) +Z_{s,k} \sin(k \omega t) ) \\
       & = l(Z_0) + l(Z_{c,1}) \cos(\omega t) + l(Z_{s,1}) \sin(\omega t) +
l(Z_{c,2}) \cos(2 \omega t) + \dots  
\end{array}
\end{equation}}
The operator $L$ is
{\footnotesize
\begin{equation}\label{opeL}\small
\hspace*{-0.2cm}
L(U)=  [ l(Z_0)^t ,  l(Z_{c,1})^t , l(Z_{s,1})^t , l(Z_{c,2})^t , l(Z_{s,2})^t , \dots , l(Z_{c,H})^t , l(Z_{s,H})^t ]^t 
\end{equation}
}
For the linear operator $m$, we have 
\begin{equation}\label{lin2}
\hspace*{-0.2cm}
\begin{array}{lll}
m(\dot Z(t))  & = & m(\sum_{k=1}^H Z_{s,k} k \omega  \cos(k \omega t) - Z_{c,k} k \omega \sin(k \omega t) ) \\
   &  = & 0 + m(Z_{s,1}) \omega \cos(\omega t) - m(Z_{c,1}) \omega \sin(\omega t)  \\
   & & + 2 m(Z_{s,2}) \omega \cos(2 \omega t) + \dots  
\end{array}
\end{equation}

The operator $M$ is
\begin{equation}\label{opeM}
\begin{array}{lll}
M(U) & = &  [ 0^t ,  m(Z_{s,1})^t , -m(Z_{c,1})^t , 2m(Z_{s,2})^t , -2m(Z_{c,2})^t , \dots \\
   &    &  \dots , Hm(Z_{s,H})^t , -Hm(Z_{c,H})^t ]^t  \;\;\;\;\;
\end{array}
\end{equation}

\subsection{Quadratic term}

With the notation $Z^{(0)}=Z_0$ and $Z^{(i)}=Z_{c,i} \cos(\omega t) + Z_{s,i} \sin(\omega t)$, the decomposition of $Z(t)$ is conveniently rewritten as
\begin{equation}\label{compactfou}
Z(t)=  \sum_{i=0}^H Z^{(i)}
\end{equation}

For the operator $q$, we have
\begin{equation}\label{devquad}
\begin{array}{ll}

q(X(t),Y(t))  &   = \sum_{i=0}^H \sum_{j=0}^H q( X^{(i)}, Y^{(j)}) \\
     & =  q( X^{(0)}, Y^{(0)}) +  \sum_{i=1}^H q( X^{(i)}, Y^{(0)})  \\
& + \sum_{j=1}^H q( X^{(0)}, Y^{(j)})   + \sum_{i=1}^H \sum_{j=1}^H q( X^{(i)}, Y^{(j)}) \\
\end{array}
\end{equation}
When $i\geq 1$ and $j\geq 1 $,  $q(X^{(i)},Y^{(j)})$ give harmonics $i+j$ and $i-j$ as follows:
\begin{equation}\label{qxy}
\begin{array}{ll}
q( X^{(i)}, Y^{(j)} ) & = \frac{1}{2}\{ q(X_{c,i},Y_{c,j}) - q(X_{s,i},Y_{s,j}) \} \cos( (i+j)\omega t) \\
    & + \frac{1}{2}\{ q(X_{c,i},Y_{s,j}) + q(X_{s,i},Y_{c,j}) \} \sin( (i+j)\omega t)\\
   & + \frac{1}{2}\{ q(X_{c,i},Y_{c,j}) + q(X_{s,i},Y_{s,j}) \} \cos( (i-j)\omega t) \\
  & + \frac{1}{2}\{ q(X_{s,i},Y_{c,j}) - q(X_{c,i},Y_{s,j}) \} \sin( (i-j)\omega t)
\end{array}
\end{equation}

By grouping the terms with the same harmonics index, and canceling any harmonics higher than index $H$, we have

\begin{equation}\label{devquadfin}
q(X(t),Y(t))  = q_0  + \sum_{k=1}^H  q_{c,k} \cos(k \omega t) +  q_{s,k} \sin(k \omega t)\\
\end{equation}
with
\begin{equation}\label{q0}
q_0  = q(X_0,Y_0) + \sum_{j=1}^{H} \frac{1}{2}\{ q(X_{c,j},Y_{c,j}) + q(X_{s,j},Y_{s,j}) \} 
\end{equation}
and when $i \geq 1$
{\scriptsize
\begin{equation}\label{devquadfin2}
\hspace*{-0.5cm}
\begin{array}{ll}
 q_{c,i}  & =  \{ q(X_{c,i},Y_0) + q(X_0,Y_{c,i}) \} +  \sum_{j=1}^{i-1}  \frac{1}{2}\{ q(X_{c,j},Y_{c,i-j}) - q(X_{s,j},Y_{s,i-j}) \}  \\
    + & \displaystyle \sum_{j=i+1}^{H}   \frac{1}{2}\{ q(X_{c,j},Y_{c,j-i}) + q(X_{s,j},Y_{s,j-i}) \}  
                     +    \frac{1}{2}\{ q(X_{c,j-i},Y_{c,j}) + q(X_{s,j-i},Y_{s,j}) \}  \\
\end{array}
\end{equation}
\begin{equation}\label{devquadfin3}
\hspace*{-0.5cm}
\begin{array}{ll}
 q_{s,i} & = \{ q(X_{s,i},Y_0) + q(X_0,Y_{s,i}) \} +  \sum_{j=1}^{i-1}  \frac{1}{2}\{ q(X_{c,j},Y_{s,i-j}) + q(X_{s,j},Y_{c,i-j}) \}  \\
    + & \displaystyle  \sum_{j=i+1}^{H}  \frac{1}{2}\{-q(X_{c,j},Y_{s,j-i}) + q(X_{s,j},Y_{c,j-i}) \}  
                     +    \frac{1}{2}\{ q(X_{c,j-i},Y_{s,j}) - q(X_{s,j-i},Y_{c,j}) \}  \\
\end{array}
\end{equation}
}

\end{document}